\documentclass[11pt]{amsart}
\usepackage{amsfonts,amssymb,amscd}

\newcommand{\bir}{\mathop{\rm {- - }\to}\nolimits}

\newcommand{\mt}[1]{\operatorname{#1}}
\newcommand{\wt}{\mt{wt}}
\newcommand{\ord}{\mt{ord}}

\newcommand{\gr}{{\mt{gr}}^0_C\omega}
\newcommand{\gro}{{\mt{gr}}^1_C\OOO}

\newcommand{\scs}{^{\mt{sc}}}

\newcommand{\bb}[1]{{\mathbb{#1}}}
\newcommand{\cc}[1]{{\mathcal{#1}}}
\newcommand{\CC}{{\bb{C}}}

\newcommand{\ZZ}{{\bb{Z}}}
\newcommand{\QQ}{{\bb{Q}}}
\newcommand{\PP}{{\bb{P}}}
\newcommand{\NN}{{\bb{N}}}
\newcommand{\OOO}{{\cc{O}}}

\newcommand{\qq}{\begin{flushright}
$\square$\end{flushright}}
\newcommand{\cyc}[1]{\bb{Z}_{#1}}
\newcommand{\3}{^{\sharp}}

\newtheorem{theorem}[subsection]{Theorem}

\newtheorem{proposition}[subsection]{Proposition}

\newtheorem{conjecture}[subsection]{Conjecture}
\newtheorem{lemma}[subsection]{Lemma}
\newtheorem{lemma1}[subsubsection]{Lemma}

\newtheorem{corollary1}[subsubsection]{Corollary}

\theoremstyle{definition}
\newtheorem{definition}[subsection]{Definition}
\newtheorem{definition1}[subsubsection]{Definition}

\newtheorem{proposition-definition}[subsection]{Proposition-Definition}
\newtheorem{proposition-definition1}[subsubsection]{Proposition-Definition}

\theoremstyle{remark}
\newtheorem{remark}[subsection]{Remark}

\newcommand{\thePOSCC}
{(\arabic{section}.\arabic{subsection}.\arabic{subsubsection})}
\newcounter{POSCC}[subsection]
\newcounter{POSCCC}[subsubsection]
\renewcommand{\thePOSCC}
{(\arabic{section}.\arabic{subsection}.\arabic{POSCC})}
\renewcommand{\thePOSCCC}
{(\arabic{section}.\arabic{subsection}.\arabic{subsubsection}.\arabic{POSCCC})}
\newcommand{\ab}{\refstepcounter{POSCC}\thePOSCC}
\newcommand{\abc}{\refstepcounter{POSCCC}\thePOSCCC}

\begin{document}
\author{Yuri G. Prokhorov}
\title{
On extremal contractions from threefolds to surfaces:\\
the case of one non-Gorenstein point and non-singular base surface
}

\email{prokhoro@mech.math.msu.su\\
prokhoro@nw.math.msu.su}

\thanks{The author was supported in part by
the Russian Foundation of Fundamental Research
\# 96-01-00820.}

\address{
Algebra Section,
Department of Mathematics,
Moscow State University,
Vorob'evy Gory, Moscow
117 234, Russia}

\date{}

\begin{abstract}
In this paper we study a local structure of
extrmal contractions
$f\colon X\to S$
from threffolds $X$ with only terminal singularities
onto a surface $S$.
If the surface  $S$  is non-singular and $X$
has a unique non-Gorenstein point on a fiber
 we prove that either the linear system
$|-K_X|$, $|-2K_X|$ or $|-3K_X|$ contains a "good"
divisor.
\end{abstract}

\maketitle

\section*{Introduction}
In this work we continue
 papers \cite{Pro}, \cite{Pro1},
devoted to study local structure of
extremal contractions
from a threffold $X$ with only terminal singularities
onto a surface $S$.

\begin{definition}\label{def}
Let $(X,C)$ be a germ of a three-dimensional complex
space   along a  compact reduced curve $C$ and let
$(S,o)$ be a germ of a
two-dimensional normal complex space.
Suppose that $X$ has  at worst terminal singularities.
Then we say that a proper morphism
$f\colon (X,C)\to (S,o)$ is a {\it Mori conic bundle}
if \par
\begin{enumerate}
\renewcommand\labelenumi{(\roman{enumi})}
\item $(f^{-1}(o))_{\mt{red}}=C$;
\item $f_*\cc{O}_X=\cc{O}_S$;
\item $-K_X$ is ample.
\end{enumerate}
 \end{definition}

As in the case of small contractions \cite{Mori-flip},
\cite{Shokurov}
 here the following conjecture is very important.

\begin{conjecture}
\label{conj}
Let $f\colon (X,C)\to (S,0)$ be a Mori conic bundle. Then for
any generic member of the linear system
 $F\in |-nK_X|$ log-canonical divisor
$K_X+\frac{1}{n}D$ is log-terminal.
\end{conjecture}

Of course, this conjecture is
interesting only for sufficiently small $n$. For instance,
 applications require it for $n=1$
(see \cite{Iskovskikh1}, \cite{Pro}).

\par
Since we are working in the analytic situation, then
applying the Minimal Model Program to  $X$ over $(S,o)$
we can achive $\rho(X/S)=1$ and then the fiber $C$ will be
irreducible and  $C\simeq\PP^1$.
In the present work we investigate the case of nonsingular base
surface $(S,o)$, irreducible fiber $C$ under an additional assumption
that $X$ contains the only one non-Gorenstein point
(in general it is known that $X$ can contain at most three
singular points\cite{Pro1} and at most two non-Gorenstein points).
It will be proved in Theorem \ref{pmain}, that
Conjecture~\ref{conj} for  $n=1$ holds, probably,
with the exception of
a finite number of exceptional cases in which
 Conjecture~\ref{conj} for $n=2$ or  $n=3$ holds.
Unfortunately, now I do not know
are these exceptional cases really possible
and if yes, then
is Conjecture~\ref{conj} for $n=1$ hold for them
(we can be sure only that in these cases
a divisor  $F\in |-K_X|$ must contain the fiber $C$).
\subsection*{Acknowledgments.}
I have been working on this problem at
Max-Planck Institut f\"ur Mathematik and
the Johns Hopkins University.
I am very grateful staffs of these institutions
for hospitality.
Different aspects of this problem were discussed
with Professors  V.I.~Iskovskikh,
V.V.~Shokurov and M.~Reid.
 I am grateful to them for help and advises.

\section{Preliminary results}
\begin{lemma}
\label{po}
Let $f:(X,C\simeq\PP^1)\to (S,0)$ be a Mori conic bundle
with non-singular base surface  $(S,o)$. Assume that $X$ is singular and
let $m_1,\dots , m_r$ be indices of
singular points on $X$ (in fact $r\le 2$).
Then $(-K_X\cdot C)=\frac{1}{m_1\cdots m_r}$.
\end{lemma}
\begin{proof}
 We claim that  the divisor class group $\mt{Cl}(X)$ is torsion-free.
Indeed if $D\in \mt{Cl}(X)$ is a torsion element, then it gives us a
cyclic cover $X'\to X$, \'etale outside $\mt{Sing}(X)$.
By taking the Stein factorization $X'\to S'\to S$ we obtain
another Mori conic bundle $X'\to S'$ and \'etale in codimension 1 cover
$S'\to S$. This contradicts smoothness of $S$.
\par
Denote by $\mt{Cl}\scs(X)$ the subgroup of Weil divisor class group
consisting of Weil divisor classes  on $X$ such that some multiplicity
of them is $\QQ$-Cartier.
If $(X,P)$ is a germ of three-dimensional terminal singularity
of index $m$, then $\mt{Cl}\scs(X,P)\simeq\cyc{m}$ (see \cite{Kawamata}).
In our case we have a natural exact sequence.
$$
\begin{CD}
0  @>>>\mt{Pic}(X)@>>>\mt{Cl}\scs(X)@>>>\oplus\mt{Cl}\scs(X,P_i)@>>>0\\
@.     @|             @.                @|\\
{} @.\ZZ           @.{}              @.\oplus\cyc{m_i}             @.{}  \\
\end{CD}
$$
where $P_i\in X$ are all the points of indices  $m_i>1$.
 Thus we have  $\mt{Cl}\scs (X)\simeq\ZZ$ and an ebbedding
   $\mt{Pic}(X)=\ZZ\hookrightarrow\mt{Cl}\scs(X)=\ZZ$ is
nothing but multiplication by $m_1\dots m_r$.
Let $D$ be the positive generator of $\mt{Cl}^{sc}(X)\simeq\ZZ$.
Then $-K_X\sim kD$ for some $k\in\NN$.
If $L=f^{-1}(s)$ is a general fiber,
then $(-K_X\cdot L)=2$. Hence $k(D\cdot L)=2$.
Since $(D\cdot L)\in\ZZ$, $k=1$ or $k=2$.
We claim that $k=1$. Assume that $k=2$. Then $(D\cdot L)=1$.
\par
There exists  a {\it standard form} of $f$
(see \cite{Beau}, \cite{Sarkisov}),
i.~e. a commutative diagram
$$
\begin{array}{ccc}
\widetilde{X}&\bir&X\\
\downarrow\lefteqn{\scriptstyle{\widetilde{f}}}&
&\downarrow\lefteqn{\scriptstyle{f}}\\
\widetilde{S}&\stackrel{\sigma}{\longrightarrow}&S\\
\end{array}
$$
where $\sigma:\widetilde{S}\to S$ is a composition of blow-ups over $o$,
$\widetilde{X}\bir X$ is a bimeromorhic map
and $\widetilde{f}:\widetilde{X}\to\widetilde{S}$ is a {\it standard} conic bundle
(in particular, $\widetilde{X}$ is non-singular and
${\widetilde{f}}^{-1}(B)$ is irreducible for any irreducible curve
$B\subset\widetilde{S}$).  Take the proper transform $\widetilde{D}$
of $D$ on $\widetilde{X}$.
For a general fiber $\widetilde{L}$ of $\widetilde{f}$ we have
$(\widetilde{D}\cdot\widetilde{L})=1$.
Since $\rho(\widetilde{X}/\widetilde{S})=1$,
$\widetilde{D}$ is $\widetilde{f}$-ample.
It gives us that each fiber of $\widetilde{f}$ is reduced and irreducible,
i.~e. the morphism $\widetilde{f}$ is smooth.
\begin{lemma1} \textup{(}see  \cite{Iskovskikh}\textup{)}
\label{l1}
Let $\pi:V\to W$ be a standard conic bundle,
let $E\subset W$ be a
$(-1)$-curve such that $\pi$ is smoothover $E$
 and let
$\sigma:W\to W'$ be the contraction of $E$.
Then there exists a diagram
$$
\begin{array}{ccc}
V&\bir&V'\\
\downarrow\lefteqn{\scriptstyle{\pi}}&
&\downarrow\lefteqn{\scriptstyle{\pi'}}\\
W&\stackrel{\sigma}{\longrightarrow}&W'\\
\end{array}
$$
where $\pi':V'\to W'$ is  a standard conic bundle and $V\bir V'$
is a bimeromorphic map that induces an isomorphism
$(V-\pi^{-1}(E))\simeq (V'-{\pi'}^{-1}(\sigma(E))$.\qq
\end{lemma1}

By Lemma \ref{l1}, there exists a standard conic bundle
$\hat{f}:\hat{X}\to S$ and a bimeromorphic map $\hat{X}\bir X$
over $S$. This map  indices an isomorphism $(\hat{X}-\hat{C})\simeq
(X-C)$, where $\hat{C}=\hat{f}^{-1}(0)$. Since $f$, $\hat{f}$ are projective
and $\rho(X/S)=\rho(\hat{X}/S)=1$,
we have $\hat{X}\simeq X$. But then $X$ is smooth, a contradiction.
\par
Thus $k=1$ and $-K_X$ is a generator of $\mt{Cl}^{sc}(X)$.
Hence $-m_1\cdots m_rK_X$ is a generator of $\mt{Pic}(X)$.
On the other hand, if we take a disc $H\subset X$ transversal to $C$,
then $(H\cdot C)=1$. Therefore $H=-m_1\cdots m_rK_X$
and $(-K_X\cdot C)=\frac{1}{m_1\cdots m_r}$.
This proves the desied assertion.
\end{proof}

\begin{lemma}
\label{irred}
Let $f:(X,C\simeq\PP^1)\to (S,o)$ be a Mori conic bundle
over a nonsingular base surface $(S,o)$ and  with
the only one non-Gorenstein point $P\in X$ of index  $m>1$.
Let $\pi\colon (X\3,P\3)\to (X,P)$ be the
$\cyc{m}$-canonical cover and $C\3:=\pi^{-1}(C)$.
Then the curve  $C\3$ is irreducible.
\end{lemma}
\begin{proof}
Indeed, assume the contrary: $C\3=C\3_1+\dots+C\3_d$.
Then the group $\cyc{m}$ permutes components of $C\3$.
Therefore  $m=dr$, $r\in\NN$.
By Lemma \ref{po} $(-K_X\cdot C)=1/m$.
On the other hand, by the projection formula
near $P$ we have
$(-K_{(X,P)}\cdot C)=\frac{1}{m}(-K_{(X\3,P\3)}\cdot C\3)=
\frac{1}{r}(-K_{(X\3,P\3)}\cdot C\3_1)$.
But $-K_{(X\3,P\3)}$ is a Cartier divisor.
Hence $(-K_{X}\cdot C)\in\frac{1}{r}\ZZ$,
a contradiction.
\end{proof}

\section{Invariants $w_P$ and $i_P$ according Mori \cite{Mori-flip}}
In this section we following Mori  \cite{Mori-flip}
recall definitions and methods of computations of
numerical invariants  $w_P$ and $i_P$ (see also \cite{Pro1}).
\subsection{}
\label{nn}
Let $(X,P)$ be threefold with only terminal singularities,
$C\subset X$ be a nonsingular curve and let $m\ge 1$ be the
index of $X$. Let $P\in X$ be an arbitrary point.
\subsubsection{}

Denote by
$\cc{I}_C$ the ideal sheaf of $C$ and $\omega_X:=\OOO_X(K_X)$.
As in  \cite{Mori-flip}, we define the
following sheaves on $C$:
$$
\begin{array}{l}
\gr:=\text{}\ \omega_X/(\cc{I}_C\omega_X),\\
\\
\gro:=\text{torsion-free part of}\ \cc{I}_C/\cc{I}^2_C.\\
\end{array}
$$
\subsubsection{}
The natural map
$$
(\omega_X\otimes\cc{O}_C)^{\otimes m}\to
\cc{O}_C(mK_X)
$$
induces an injection
$$
\beta\colon  (\gr)^{\otimes m}\to \cc{O}_C(mK_X).
$$
Denote
$$
w_P:=(\mt{length}_P\mt{Coker}\beta)/m.
$$
\subsubsection{}
We also have the natural  map
$$
\begin{CD}
\alpha\colon \wedge^2(\gro)\otimes\omega_C
@>>>
\omega_X\otimes\cc{O}_C
\to
\gr,\\
x\wedge y\times zdt@>>> zdx\wedge dy\wedge dt\\
\end{CD}
$$
Let
$$
i_P:=\mt{length}_P\mt{Coker}(\alpha ).
$$
\begin{lemma1}
\textup{(}\cite[2.15]{Mori-flip}\textup{)}
\label{non}
Point  $(X,P)$ is singular iff  $i_P\ge 1$.\qq
\end{lemma1}

\begin{proposition}
\textup{(} see \cite{Mori-flip}, \cite{Pro1}\textup{)}
\label{grw}
Let $f\colon (X,C\simeq\PP^1)\to (S,0)$ be a Mori conic bundle.
Then
$$
(-K_X\cdot C)+\sum w_P+\sum i_P\le 4.
\eqno\square
$$
\end{proposition}

As in  \cite{Mori-flip}
we need the following construction to
compute local invariants $w_P$ and  $i_P$.
\subsection{}
\label{nnn}
In notations of  \ref{nn}
consider the canonical $\cyc{m}$-cover
$\pi\colon (X\3,P\3)\to (X,P)$
and denote $C\3:=\pi^{-1}(C)$.
Later throughout in this paper we will assume that
the curve $C\3$ is irreducible
(this holds in our case by Lemma \ref{irred}).
There exists an $\cyc{m}$-equivariant embedding $X\3\subset\CC^4$.
Let $\phi=0$ be the equation of $X\3$ in $\CC^4$.
\begin{definition1}
Fix some character
$\chi\colon \cyc{m}\to\CC^*$ with trivial kernel.
Then we define {\it weight} of semi-invariant $x$
as an integer $\wt(x)$ defined $\mod m$
 such that
$$
\xi (x)=\chi(\xi)^{\wt(x)}\cdot x
\qquad \text{для всех}
\xi\in\cyc{m}.
$$
\end{definition1}
\subsubsection{}
\label{cl.term}
By the classification of terminal singularities \cite{Mori-term},
if  $(x_1,x_2,x_3,x_4)$ is a semi-invariant coordinate system
in $\CC^4$, then up to permutation one of the following cases
below holds
\begin{enumerate}
\renewcommand\labelenumi{(\roman{enumi})}
\item
(the main series)
$\wt(x_1)+\wt(x_2)\equiv\wt(x_4)\equiv\wt(\phi)\equiv 0\bmod m$,
$\wt(x_i)$ is prime to $m$ for $i=1,2,3$;
\item
(the exceptional series)
$m=4$,
$\wt(x_1)+\wt(x_2)\equiv 0\bmod 4$,
$\wt(x_4)\equiv\wt(\phi)\equiv 2\bmod 4$ and
$\wt(x_i)$ is odd for $i=1,2,3$.
\end{enumerate}
There is the full list of normal forms of
equations $\phi$ for all terminal singularities,
usually we distinguish the following cases
by type of generic section $X\3\cap\CC^3\ni 0$:
$cA/m$, $cAx/2$, $cD/2$, $cD/3$, $cE/2$
in the main series и and the only case
$cAx/4$ in the exceptional series.
\begin{definition1}
For $z\in\OOO_{X\3}$ define its {\it order} $\ord(z)$
as the order of vanishing of $z$ on the normalization of $C\3$.
All the values $\ord(z)$ form a semigroup which denoted by
$\ord(C\3)$. Obviously $\ord(x_1),\dots,\ord(x_4)$
generate  $\ord(C\3)$.
\end{definition1}

\begin{proposition-definition}
\label{notations}
\label{nor}
By \cite{Mori-flip} one can
choose a coordinate system
 $(x_1, x_2, x_3, x_4) $  in
$\CC^4$  and a character $\chi\colon \cyc{m}\to\CC^*$
such that the following conditions hold.
\begin{enumerate}
\renewcommand\labelenumi{(\roman{enumi})}
\item
Locally near $P\3$ the curve  $C\3$ is the image of the map
$$
\CC\to \CC^4,\qquad
t\longrightarrow
(t^{\ord(x_1)},t^{\ord(x_2)},t^{\ord(x_3)},t^{\ord(x_4)}),
$$
where $t$ is a semi-invariant with  $\wt(t)\equiv 1\bmod m$.
In particular, all the coordinates $x_1, x_2, x_3, x_4$
are semi-invariants and
$\ord(x_i)\equiv\wt(x_i)\bmod m$ for all $i=1,2,3,4$;
\item\label{normal}
{\it Normalizedness property.}
$\ord(x_i)<\infty$ for all  $i$ and
there is no semi-invariants $y$
such that  $\wt(y)\equiv\wt(x_i)\bmod m$ and $\ord(y)<\ord(x_i)$.
\item
There exists an invariant function $z$ on $X\3$ such that
$\ord(z)=m$.
\end{enumerate}
Such a coordinate system is said to be {\it normalized}.
\end{proposition-definition}
\begin{remark}
\label{notations1}
By permutation of coordinates we can acheve that
weights $\wt(x_i)$ are sutisfied (i) or (ii) of
\ref{cl.term}.
In particular, if $(X,P)$ is from the matn series, then
by normalizedness and (ii) \ref{nor} we have
$\ord(x_4)=m$, $\wt(x_4)\equiv 0\bmod m$.
In the case of singularity of exceptional series
since $4\in\ord (C\3)$, we have $\ord (x_4)=2$.
Note that we still may  permute $x_1$, $x_2$.
\end{remark}

\begin{proposition}
\textup{(} see \cite[2.10]{Mori-flip}\textup{)}
\label{computation-w}
Under
  \textup{\ref{notations}}, \textup{\ref{notations1}}
 one has
$$
w_P=\min \{ \ord(\psi)\quad |\quad  \psi\in\OOO_{X\3,P\3},\quad
\wt(\psi)\equiv-\wt(x_3)\bmod m\}.
$$
\qq
\end{proposition}

\begin{proposition}
 \textup{(} see \cite[0.4.14.2]{Mori-flip}\textup{)}
\label{computation-k}
Let  $F\in |-K_{(X,P)}|$ be a general member.
the notations and conditions  of
\textup{\ref{notations}}, \textup{\ref{notations1}}
one has
$(F\cdot C)_P=\ord(x_3)/ m$.
\qq
\end{proposition}
As an immediate consequence of Lemma \ref{po} we have.
\begin{corollary1}
\label{a3p}
Let $f\colon (X,C\simeq\PP^1)\to (S,0)$ be
a Mori conic bundle over a non-singular base surface
 $(S,o)$ with a unique non-Gorenstein point
 $P\in X$ of index  $m$.
Then in notations above we have
$\ord(x_3)\equiv 1\bmod m$.
\qq
\end{corollary1}

\begin{proposition}
\textup{(} see \cite[2.12]{Mori-flip}\textup{)}
\label{computation-i}
In the notations and conditions  of
\textup{\ref{notations}}, \textup{\ref{notations1}}
one has
$$
i_Pm=m-\ord(x_4)-m w_P+\min_{\phi_1,\phi_2\in
J\3_0}[\phi,\phi_1,\phi_2],
$$
where $J\3_0$  is the invariant part of the ideal sheaf
$J\3$ of curve $C\3$ in $\CC^4$, $\phi$ is an equation of
$X\3$ in  $\CC^4$ and
$$
[\psi_1,\psi_2,\psi_3]:=
\ord\partial (\psi_1,\psi_2,\psi_3)/\partial (x_1,x_2,x_3).
$$
is the Jacobian determinant.\qq
\end{proposition}
An invariant monomial $\psi$ in $x_1, x_2, x_3$ is said
to be {\it simple}
if it cannot be presented as a product of two non-constant
invariant monomials.
By semi-additivity of $[\ ,\ ,\ ]$ and
because $x_4\phi\in J\3_0$, we have
\begin{corollary1}
\textup{(} see. \cite[ proof of 2.15]{Mori-flip}\textup{)}
\label{formula}
\label{computation-i1}
In the notations and conditions  of
\textup{\ref{notations}}, \textup{\ref{notations1}}
for some simple different invariant monomials
$\psi_1$, $\psi_2$, $\psi_3$
 in  $x_1, x_2, x_3$ we have
$$
\begin{array}{l}
m i_P\ge(\ord(\psi_1)-\ord(x_1))+(\ord(\psi_2)-\ord(x_2))+\\
\qquad (\ord(\psi_3)-\ord(x_3)-m w_P).
\end{array}
\label{(*)}
\leqno\abc
$$
 Moreover up to permutations $\psi_1$, $\psi_2$, $\psi_3$
we may assume that $\psi_i=x_i\nu_i$, where $\nu_i$, $i=1,2,3$ also
are  monomials. In particular,
all three terms in the formula are non-negative and
$$
m i_P\ge\ord(\nu_1)+\ord(\nu_2)+\ord(\nu_3)-m w_P.
\label{(**)}
\leqno\abc
$$
\qq
\end{corollary1}

\section{The main result}
 In this section we prove the following
\begin{theorem}
\label{pmain}
Let $f:(X,C\simeq\PP^1)\to (S,0)$ be a Mori conic bundle with
irreducible central fiber. Suppose that $X$ contains only one non-Gorenstein
point $P$ and $(S,0)$ is non-singular. Let $m$ be the index of $P$.
Then
general member of $|-K_X|$ does not contain $C$ and has
DuVal singularity at $P$, except possibly one of the following cases
in which a general member  $F\in |-2K_X|$ does not contain $C$
\textup{(}cases \textup{\ref{2-cAx/4}}, \textup{\ref{2-IC}}
 and \textup{\ref{2-cA/m})} or
a general member  $F\in |-2K_X|$ does not contain $C$
\textup{(}cases
\textup{\ref{3-7}, \ref{3-cyc1}, \ref{3-diff1},
\ref{3-cyc2}, \ref{3-diff2} }
and \textup{\ref{3-diff3})}.
\subsubsection{}
\label{2-cAx/4}
The point $(X,P)$ is of type $cAx/4$ and near  $P$ we have
$$
(X\supset C)\simeq (\{y_1^2-y_4^3+y_2\varphi_1+y_3\varphi_2\}
\supset \{y_2=y_3=y_1^2-y_4^3=0\}/\cyc{4}(3,1,1,2),
$$
where
$\wt(\varphi_1)\equiv\wt(\varphi_2)\equiv 1\bmod 4$
и $y_2\varphi_1+y_3\varphi_2\not\in (y)^3$
\textup{(} a point of type $IIB$ in notations of~\cite{Mori-flip}\textup{)}.
\subsubsection{}
\label{2-IC}
$m$ is odd,  $m \ge 5$ and near  $P$ we have
$$
(X\supset C)\simeq
(\CC^3\supset \{y_3=y_2^2-y_1^{m-2}=0\})/\cyc{m}(2,-2,1)
$$
\textup{(}a point of type  $IC$ in notations \cite{Mori-flip}\textup{)}.
\subsubsection{}
\label{2-cA/m}
The point $(X,P)$ is of type  $cA/m$, $m$ is odd and
$$
(X,P)\simeq\{y_1y_2+y_2\varphi_1+y_3\varphi_2+(y_4^2-y_1^{m})\varphi_3=0\}
/\cyc{m}(2,-2,1,0),
$$
where \qquad
$
C=(\{y_2=y_3=y_4^2-y_1^{m}=0\}/\cyc{m}.
$
\subsubsection{}
\label{3-7}
$m=7$,\quad
$(X,P)\simeq (\CC^3,0)/\cyc{7}(3,-3,1)$, where
$C=\{ y_3=y_1^4-y_2^3=0\}/\cyc{7}$.
\subsubsection{}
\label{3-cyc1}
$m\ge 10$,\quad $m\equiv 1\bmod 3$ and
$(X,P)\simeq (\CC^3,0)/\cyc{m}(3,-3,1)$, where
$
C=\{y_2^{3}-y_1^{m-3}=y_3^3-y_1^{m+1}=0\}/\cyc{m}.
$
\subsubsection{}
\label{3-diff1}
 $m\equiv 1\bmod 3$,
the point $(X,P)$ is of type $cA/m$ and
\begin{multline*}
(X,P)\simeq
\{x_1x_2-x_4^2+(x_2-x_3x_1^{(m-4)/3})\varphi_2+
(x_1^{m+1}-x_3^3)\varphi_3+\\
(x_1^m-x_4^3)\varphi_4+
(x_3^m-x_4^{m+1})\varphi_5 =0\}/\cyc{m}(3,-3,1,0),
\end{multline*}
where
$$
C=\{ x_1x_2-x_4^2=x_2-x_3x_1^{(m-4)/3}=x_1^{m+1}-x_3^3=
x_1^m-x_4^3=x_3^m-x_4^{m+1}=0\}/\cyc{m}.
$$
\subsubsection{}
\label{3-cyc2}
$m\ge 10$,  $m\equiv 1\bmod 3$ and near  $P$ we have
$$
(X\supset C)\simeq
(\CC^3\supset\{y_2^3-y_1^{m-3}=y_3=0\})/\cyc{m}(3,-3,1).
$$
\subsubsection{}
\label{3-diff2}
 $m\equiv 1\bmod 3$,
$m\ge 7$ and
$$
(X,P)\simeq
\{y_1y_2-y_4^2+y_3\varphi_3+(y_1^m-y_4^3)\varphi_4+
(y_2^m-y_4^{2m-3})\varphi_5=0\}/\cyc{m}(3,-3,1,0),
$$
$$
C=\{y_1y_2-y_4^2=y_3=y_1^m-y_4^3=y_2^m-y_4^{2m-3}=0\}/\cyc{m}.
$$
\subsubsection{}
\label{3-diff3}
 $m\equiv 1\bmod 3$,
the point $(X,P)$ is of type $cA/m$ and
$$
(X,P)\simeq\{y_1y_2+y_2\varphi_2+y_3\varphi_3+
(y_1^m-y_4^3)\varphi_4\}/\cyc{m}(3,-3,1,0),
$$
where \qquad
$
C=\{y_2=y_3=y_1^m-y_4^3=0\}/\cyc{m}.
$
\par
Moreover in cases
\textup{\ref{2-cAx/4}, \ref{3-7}, \ref{3-cyc1},
\ref{3-diff1}, \ref{3-cyc2}, \ref{3-diff2} }
and \textup{\ref{3-diff3}} $X$ has no other (Gorenstein) singular
points and in cases \textup{\ref{2-IC},\ref{2-cA/m}}
 $X$ may has at most one Gorenstein singular point.
\end{theorem}

\subsection{}\label{asu}
The reminder of this paper is devoted to proof of this theorem.
Throughout this section let  $f:(X,C)\to (S,o)$ be
Mori conic bundle such as in Theorem~\ref{pmain}, $P$ be
a (unique) non-Gorenstein point on $X$ and $m$ be its index.
We also will use all the notations of
\ref{nn}, \ref{nnn}, \ref{notations}.
Set $a_i:=\ord(x_i)$.
To prove that a general
member of $|-K_X|$ does not contain $C$,
by Proposition \ref{computation-k},
it is sufficient to
show that $a_3<m$ (in fact, then we have
$a_3=1$ by Corollary \ref{a3p}).   So we asusume that
$a_3>m$ and then from the normalizedness we have
$1\not\in\ord(C\3)$. Our main tool will be inequality \ref{grw}.
\par
First, we consider the case of singularity of exceptional series.

\begin{lemma}
\label{uuu}
\textup{(} cf. \cite[(4.2), (4.4)]{Mori-flip}\textup{)}
Let in assumptions of \textup{\ref{asu}} $(X,P)$ be a
terminal point of type $cAx/4$.
 Then
 $2\in \ord(C\3)$ and near $P$
\begin{multline*}
(X\supset C)\simeq \\
(\{y_1^2-y_4^3+y_2\varphi_1+y_3\varphi_2=0\}
\supset \{y_2=y_3=y_1^2-y_4^3=0\})/\cyc{4}(3,1,1,2),
\end{multline*}
 where
$\varphi_1$, $\varphi_2$ are semi-invariants with
$\wt(\varphi_1)\equiv\wt(\varphi_2)\equiv 1\bmod 4$.
Moreover $i_P=3$, $w_P=3/4$ in this case.
\end{lemma}

\begin{proof}
By our assumptions $1\not\in\ord(C\3)$.
From normalizedness  and because $(a_i,2)=1$, $i=1,2,3$,
we may assume that
$a_2=a_3$ up to permutatation $\{1, 2\}$.
We also have $a_3>4$, $a_3\equiv 1\bmod 4$.
Hence $a_1\equiv 3\bmod m$.
\par
First, we assume that  $a_1=3$.
Since $a_4=2$ by \ref{notations1},
we have $5\in\ord(C\3)$, hence $a_2=a_3=5$.
Thus $a_1=3$, $a_2=a_3=5$, $a_4=2$
 and $C\3$ is the image of $t\longrightarrow (t^3, t^5, t^5, t^2)$.
Changing coordinates by $y_2=x_2-x_1x_4$, $y_3=x_3-x_2$, $y_1=x_1$, $y_4=x_4$
 we obtain
that $C\3$ in a new (non-normalized) coordinate system can
 be given by $y_2=y_3=y_1^2-y_4^3=0$.
Therefore $X\3=\{(y_1^2-y_4^3)\varphi_0+y_2\varphi_1+y_3\varphi_2=0\}$.
But up to permutation $\{2,3\}$
an equation of $X\3$ must contain term $y_1^2+y_2^2$ (see \cite{Mori-term}),
Hence $\varphi_0$ is a unit and we may take $\varphi_0=1$.
Finally $i_P=3$ follows from \ref{computation-i1}
Thus we have the situation in Lemma.
\par
Consider now the case $a_1\ge 7$. Since $a_1, a_2, a_3>4$,
Proposition \ref{computation-w} gives us
$w_P>1$. Therefore $i_P\le 2$.
By \ref{formula}, $\nu_1+\nu_2\le 8$, where
$\nu_1$, $\nu_2$ are semi-invariants with
 $\ord(\nu_1)\equiv a_2\bmod 4$ and
$\ord(\nu_2)\equiv a_1\bmod 4$. From normalizedness we have $a_1+a_2\le 8$,
a contradiction.
\end{proof}

Therefore if $(X,P)$ is a singularity from the exceptional series,
then we get the case \ref{2-cAx/4}.
Consider now the case of singularities from the main series.

\begin{lemma}
\label{23}
In assumptions of \textup{\ref{asu}},
one has $2\in\ord(C\3)$ or $3\in\ord(C\3)$.
\end{lemma}

\subsection{Proof}\label{label}
Assume that $1,2,3\not\in\ord(C\3)$. Then $m\ge 4$ and by Lemma
above $(X,P)$ is from the  main series.
Recall also that $a_3\equiv 1\bmod m$ and $a_3\ge m+1$.
We have by Corollary \ref{formula}
$$
\ord(\nu_1)+\ord(\nu_2)+\ord(\nu_3)-mw_P\le 3m,
\label{(***)}
\leqno\ab
$$
where $\nu_i$, $i=1,2,3$ are
monomials such that $\psi_i=\nu_ix_i$.
If $\psi_1=x_1^m$, then $\ord(\nu_1)=(m-1)a_1< 3m$.
Whence $a_1< 3m/(m-1)\le 4$, a contradiction.
Thus later we will assume that $\psi_1\ne x_1^m$ and
similarly $\psi_2\ne x_2^m$.

By normalizedness,
$$
a_1+a_2+\ord(\nu_3)-mw_P\le \ord(\nu_1)+\ord(\nu_2)+\ord(\nu_3)-mw_P\le 3m.
$$
Since $a_1+a_2\equiv 0\bmod m$, $a_1+a_2=m$, $2m$ or $3m$.
Consider these cases.
\subsection{Case $a_1+a_2=3m$}
\label{ou}
Then $\ord(\nu_3)=mw_P$,
$\ord(\nu_1)+\ord(\nu_2)=3m$ and we have an equality in
\ref{(***)}. Therefore $i_P=3$ and $w_P<1$ (see \ref{grw}).
Thus $\ord(\psi_3)=\ord(\nu_3)+a_3<m+a_3$.
Up to permutation of $\{1, 2\}$ we may suppose that $a_1<a_2$, so
$a_2>m$ and by our assumption $a_3>m$.
It gives us $\nu_3=x_1^{\gamma}$ and then $\psi_3=x_1^{\gamma}x_3$,
where $a_1\gamma<m$, $a_1\gamma +a_3\equiv 0\bmod m$.
Futher $\ord(\nu_2)=a_1<m$, hence $\nu_2=x_1$ and
$\psi_2=x_1x_2$. Since
$$
\ord(\psi_1)+\ord(\psi_2)=\ord(\nu_1)+\ord(\nu_2)+a_1+a_2=6m,
$$
$\ord(\psi_1)=3m$.
On the other hand $a_2=3m-a_1>2m$.
Therefore
$$
\psi_1\in\{x_1^{\alpha}x_3,\quad  x_1^{\alpha}x_3^2\}.
$$
But if $\psi_1=x_1^{\alpha}x_3^2$, then $a_1\alpha+2a_3=3m$, $a_3<2m$.
So $a_1\gamma+a_3=2m$ and $a_1(2\gamma-\alpha)=m$.
Since $(a_1, m)=1$, we have $a_1=1$, a contradiction with our assumptiion
$1\not\in\ord (C\3)$.
Therefore $\psi_1=x_1^{\alpha}x_3$, where $a_1\alpha+a_3=3m$, $a_3<3m$.
As above we have $a_1\gamma+a_3=2m$ or $3m$ and $a_1(\alpha-\gamma)=m$ or $0$.
Thus either $a_1=1$ or $\alpha=\gamma$ (and then $\psi_1=\psi_3$),
a contradiction.

\subsection{Case $a_1+a_2=2m$}
Up to permutation we may suppose $a_1<m$ and $a_2>m$.
In this case $\ord(\nu_1)+\ord(\nu_2)=2m$ or $3m$.
\subsubsection{Subcase $\ord(\nu_1)+\ord(\nu_2)=2m$}
Then  $\ord(\psi_1)=\ord(\psi_2)=2m$.
 As in the previous case
$\psi_2=x_1x_2$, $\psi_1=x_1^{\alpha}x_3$, where
$\alpha a_1+a_3=2m$, $a_3<2m$.
Hence $\wt(x_1^{\alpha})\equiv-\wt(x_3)\bmod m$ and
$\ord(x_1^{\alpha})<m$. By taking in  \ref{computation-w}
$\psi=x_1^\alpha$, we obtain $w_P<1$.
So $\ord(\psi_3)\le a_3+m+mw_P\le 3m$.
Whence
$$
\psi_3\in\{ x_2x_3,\quad   x_1^{\gamma}x_3,\quad    x_1^{\gamma}x_3^2\}.
$$
Similar to \ref{ou}
$\psi_3\ne x_1^{\gamma}x_3,\quad    x_1^{\gamma}x_3^2$.
But if $\psi_3=x_2x_3$, then $\wt(x_3)\equiv\wt(x_1)$ and
by normalizedness $a_3=a_1<m$, a contradiction.
\subsubsection{Subcase $\ord(\nu_1)+\ord(\nu_2)=3m$, i.~e.
$\ord(\psi_1)+\ord(\psi_2)=5m$}
Then by \ref{(***)} it is easy to see that
$\ord(\psi_3)=a_3+mw_P$, $i_P=3$, $w_P<1$ and  $\ord(\psi_3)<a_3+m$.
As in \ref{ou} we have $\psi_3=x_1^{\gamma}x_3$, where
$a_1\gamma+a_3\equiv 0\bmod m$ and $a_1\gamma<m$.
Since $\ord(\psi_2)\le 3m$ and $a_2, a_3>m$, we have
$\psi_2\in\{ x_1x_2,\quad   x_2x_3\}$. But if $\psi_2=x_2x_3$, then
$\wt(x_3)\equiv\wt(x_1)$, hence $a_3=a_1<m$, a contradiction.
Therefore $\psi_2=x_1x_2$ and  $\ord(\psi_1)=3m$.
Then $\psi_1\in\{x_1^{\alpha}x_3,\quad   x_1^{\alpha}x_3^2\}$.
Both possibilities give us  contradictions as in \ref{ou}.

\subsection{Case $a_1+a_2=m$}
Then $\ord(\psi_1)+\ord(\psi_2)\le 4m$.
It is clear that $\ord(\psi_1)+\ord(\psi_2)\ge 3m$.
\subsubsection{Subcase $\ord(\psi_1)+\ord(\psi_2)=3m$}
Up to permutation we have $\ord(\psi_1)=m$, $\ord(\psi_2)=2m$.
Hence $\psi_1=x_1x_2$, $\psi_2=x_2^{\beta}x_3$, where
$$
a_2\beta+a_3=2m\qquad \beta\ge 2,
\label{(****)}
\leqno\abc
$$
(the last follows from normalizedness).
 It gives us $a_3<2m$ and
$\wt(x_2^{\beta})\equiv -\wt(x_3)\bmod m$,
$a_2\beta<m$, and then by Proposition \ref{computation-w} $w_P<1$.
Whence by \ref{(***)} $\ord(\psi_3)\le mw_P+a_3\le 3m$.
We get the following possibilities
$$
\psi_3\in\{ x_1^{\gamma}x_3,\quad   x_1^{\gamma}x_3^{2}, \quad
x_2^{\gamma}x_3,\quad   x_2^{\gamma}x_3^{2}\}.
$$
If $\psi_3=x_2^{\gamma}x_3$, then $a_2\gamma+a_3=2m$ or $3m$.
By using \ref{(****)} we derive $a_2(\gamma-\beta)=0$ or $m$.
It follows from $(a_2, m)=1$
that $\gamma=\beta$ or $a_2=1$, a contradiction.
\par
If $\psi_3=x_2^{\gamma}x_3^2$, then
 $a_2\gamma+2a_3=3m$,
$a_2\beta +a_3=2m$.
As above we have $(2\beta-\gamma)a_2=m$, but then
$2m=a_2\frac{\gamma+m}{2}+a_3>a_2\frac{\gamma+m}{2}+m$
and  $a_2<2$,
a contradiction.
\par
If $\psi_3=x_1^{\gamma}x_3$, then $a_1\gamma+a_3=2m$ or $3m$ and
$a_2\beta +a_3=2m$. From normalizedness we have $\beta, \gamma\ge 2$.
Then $a_1(\gamma-2)+a_2(\beta-2)+2a_3=3m$.
It gives us $\gamma=2$ or $\beta=2$.
Since $a_3\equiv 1\bmod m$ and $a_3<2m$, $a_3=m+1$.
Whence $a_1\gamma=2m-1$, $a_2\beta=m-1$ and $\gamma\ne 2$.
Therefore $\beta=2$, $a_2=(m-1)/2$, $a_1=m-a_2=(m+1)/2$.
Thus $(m+1)\gamma/2=2m-1$, $m(\gamma-4)+\gamma+2=0$.
We obtain $\gamma=3$, $m=5$, $a_2=2$, a contradiction with
our assumption $2\not\in\ord (C\3)$.
\par
Finally, if $\psi_3=x_1^{\gamma}x_3^2$, then $a_1\gamma+2a_3=3m$ and
$a_2\beta +a_3=2m$. So $a_1(\gamma-1)+a_2(\beta-1)+3a_3=4m$.
  Therefore $\gamma=1$.
Since $a_3=m+1$, $a_1=m=2$ and $a_2=2$, a contradiction.
\subsection{}
Thus later we may assume that
$\ord(\psi_1)+\ord(\psi_2)=4m$.
Then by \ref{(*)} $i_P=3$, $w_P<1$, $\ord(\psi_3)=a_3+mw_P<a_3+m$.
Whence we have only two possibilities
 $\psi_3\in \{x_1^\gamma x_3,\quad
x_2^\gamma x_3\}$, where $\gamma\ge 2$ (by normalizedness).
Permuting  $x_1, x_2$, if neseccary, we obtain
$\psi_3=x_1^\gamma x_3$, where
$$
a_1\gamma=m-1, \qquad \gamma\ge 2.
\label{posl}
\leqno\ab
$$

\subsubsection{Subcase $\ord(\psi_1)=m$, $\ord(\psi_2)=3m$}
Then $\psi_1=x_1x_2$ and $\psi_2\in\{ x_2^\beta x_3^2,\quad
x_2^\beta x_3\}$.
\par
If $\psi_2=x_2^\beta x_3^2$, then $a_2\beta+2a_3=3m$,
$a_3<2m$ and $a_3=m+1$.
Thus $a_2\beta=m-2$. On the  other hand,  by \ref{posl}
we have $a_1(\gamma-1)+a_2(\beta-1)=m-3$, so $\beta=1$, $a_2=m-2$
and $a_1=2$,  a contradiction.
\par
Therefore $\psi_2=x_2^\beta x_3$, $a_2\beta+a_3=3m$, $a_3<3m$,
and $a_3=m+1$ or  $2m+1$.
Then $a_2\beta=3m-a_3=2m-1$ or $m-1$.
By normalizedness $\beta>1$.
By \ref{posl}, $0\le a_1(\gamma-2)+a_2(\beta-2)=2m-a_3-1$ and it
 gives us
$a_3=m+1$ and  $\gamma=2$ or $\beta=2$.
But if  $\gamma=2$, then $a_1=(m-1)/2$ and $a_2=(m+1)/2$.
This contradicts to  $a_2\beta=(m+1)\beta/2=2m-1$.
Thus  $\beta=2$ and then $2a_2=2m-1$, again a contradiction.

\subsubsection{Subcase $\ord(\psi_1)=3m$, $\ord(\psi_2)=m$}
Then $\psi_2=x_1x_2$ and $\psi_1\in\{ x_1^\alpha x_3^2,\quad
x_1^\alpha x_3\}$. But if $\psi_1= x_1^\alpha x_3$,
then $\psi_1$ is divisible by $\psi_3$.
This is impossible because monomials $\psi_1,\psi_3$ are simple.
Therefore $\psi_1= x_1^\alpha x_3^2$ and then
$a_3=m+1$, $a_1\alpha=m-2$, this contradicts \ref{posl}.

\subsubsection{Subcase $\ord(\psi_1)=\ord(\psi_2)=2m$}
Then
$$
\psi_1, \psi_2\in\{x_1x_2, x_1^{\alpha}x_3, x_2^{\beta}x_3\}.
$$
It gives us $a_3<2m$ and $\ord(\psi_3)=2m$. Hence we also have
$$
\psi_3\in\{x_1x_2, x_1^{\alpha}x_3, x_2^{\beta}x_3\}.
$$
Up to permutation we have
$\psi_1=x_1x_2$, $\psi_3=x_1^{\alpha}x_3$ and
$\psi_2=x_2^{\beta}x_3$, where
$a_1\alpha+a_3=2m$ and $a_2\beta+a_3=2m$. As above we get
$\alpha=1$ or $\beta=1$, a contradiction.
This proves Lemma \ref{23}.
\qq
\par

Now Theorem~\ref{pmain} will follow from the following two Lemmas.

\begin{lemma}
\label{cb}
Let in assumptions of \textup{\ref{asu}}
$2\in\ord(C\3)$, $1\not\in\ord(C\3)$
and $(X,P)$ is a point from the main series.
Then $m$ is odd, up to permutation of $x_1$, $x_2$, $a_1=2$  and
we have one of the following two cases in the table.
$$
\begin{array}{c||c|c|c|c|c}
\mt{}&a_1 &a_2 &a_3&m&i_P\\
\hline\\
\ab\label{2-1)}         &2   &m-2 &m+1&\mt{odd},\quad m\ge 5&2\\
\ab\label{2-2)}        &2   &2m-2&m+1&\mt{odd}             &\ge 2\\
\end{array}
$$
Moreover in the case \textup{\ref{2-1)}} $(X,P)$ is a cyclic quotient singularity,
and in the case \textup{\ref{2-2)}} $(X,P)$ is of type $cA/m$.
\end{lemma}
\begin{proof}
We claim that $a_1=2$. Indeed, if otrherwise,
$m=a_4=2$ by assumption $a_3>m$. But then $a_1=a_2=a_3\ge 3$.
Hence $w_P>1$, $i_P\le 2$. By Corollary~\ref{formula} we have
$\ord(\nu_1)+\ord(\nu_2)\le 2m=4$, where
$\ord(\nu_1)\equiv a_2$, $\ord(\nu_2)\equiv a_1\bmod 2$.
From normalizedness we have $2a_1=a_1+a_2\le 4$, a contradiction.
\par
Therefore $a_1=2$ and $m$ is odd.
By Corollary~\ref{a3p}, $a_3\equiv 1\bmod m$.
Then $\ord(C\3)$ is generated by
$2$ and the smallest odd $c\in\ord(C\3)$. Two possibilities
$c=a_2$ and $c=m$ give us two cases in the table.
Finally in the case \ref{2-1)} $C\3$ is the image of
$t\longrightarrow (t^{2},t^{m-2},t^{m+1},t^{m})$.
If we change coordinates by $x_3'=x_3-x_1^{(m+1)/2}$, $x_4'=x_4-x_1x_2$,
we obtain that $C\3=\{x_3=x_4=x_2^2-x_1^{m-2}=0\}$.
Whence an equation of $X\3$ is
$\phi=x_3\varphi_1+x_4\varphi_2+(x_2^2-x_1^{m-2})\varphi_3=0$.
If $X\3$ is singular, then this equation must contain term $x_1x_2$,
because $m\ge 5$ and $\wt(x_3)\not\equiv\wt(x_1)$,
$\wt(x_3)\not\equiv\wt(x_2)$. This is a contradiction.
Similarly in the case \ref{2-2)}
$\phi=x_2\varphi_1+x_3\varphi_2+(x_4^2-x_1^{m})\varphi_3=0$.
If $(X,P)$ is not of type $cA/n$, then it is of type $cD/3$
and $\phi$ must contain the term $x^2_4$, a contradiction.
Values can be cumputed from \ref{computation-i1}.
\end{proof}

\begin{lemma}
\label{cv}
Let in assumptions of \textup{\ref{asu} }
 $3\in\ord(C\3)$ and $2\not\in\ord(C\3)$.
Then
\subsubsection{}
\label{t1}
$i_P=3$,\quad $w_P=(m-1)/m$.

\subsubsection{}
\label{t2}
$m\equiv 1\bmod 3$, the point $(X,P)$ is of type $cA/m$
and up to permutation of $\{1, 2\}$
we have $a_1=3$.

\subsubsection{}
\label{t3}
Only one of the following cases in the table below holds.
$$
\begin{array}{c||c|c|c|l}
\mt{}&a_1&a_2&a_3&m\\
\hline\\
\abc \label{1)}     &3&m-3  &m+1   &m\equiv 1\bmod 3,\quad m\ge 7    \\
\abc \label{2)}     &3&2m-3 &m+1  &m\equiv 1\bmod 3                 \\
\abc \label{3)}     &3&m-3  &2m+1  &m\equiv 1\bmod 3,\quad  m\ge 10  \\
\abc \label{4)}     &3&2m-3 &2m+1  &m\equiv 1\bmod 3,\quad m\ge 7    \\
\abc \label{5)}     &3&3m-3 &2m+1  &m\equiv 1\bmod 3                 \\
\end{array}
$$
Moreover in cases \textup{\ref{1)},\ref{3)}}
 $(X,P)$ is a cyclic quotient.
In case  \textup{\ref{5)}} $(X,P)\simeq\{y_1y_2+y_2\varphi_2+y_3\varphi_3+
(y_1^m-y_4^3)\varphi_4\}/\cyc{m}$.
\end{lemma}
\begin{proof}
Since $2\not\in\ord(C\3)$, $P$ is from  the main series (see
\ref{uuu}).
Similar to proof of Lemma \ref{cb} it is easy to
see that $a_1=3$. Suppose that $i_P\le 2$.
Then by Corollary \ref{formula}, for some simple invariants
$\psi_i$ we have
$$
2m\ge (\ord(\psi_1)-a_1)+(\ord(\psi_2)-a_2)+(\ord(\psi_3)-a_3-mw_P).
$$
Similar to \ref{label}, $a_1+a_2\le 2m$.

\subsection{Case $a_1+a_2=2m$}
Then $\ord(\psi_1)+\ord(\psi_2)=4m$.
Up to permutation $\{1,2\}$ we have
$\psi_1=x_1^\alpha x_3$ and $\psi_2=x_1x_2$, where $a_1\alpha+a_3=2m$.
Therefore $a_3=m+1$.
By Lemma \ref{computation-w}, $mw_P=\ord(x_1^\alpha)<m$.
Whence $\ord(\psi_3)=a_3+mw_P=2m$.
But then $\psi_3=\psi_1$, a contradiction.

\subsection{Case  $a_1+a_2=m$}
Then again $\ord(\psi_1)+\ord(\psi_2)=3m$,
$\psi_1=x_1^\alpha x_3$ and $\psi_2=x_1x_2$, where $a_1\alpha+a_3=2m$.
Therefore $a_3=m+1$.
Again by Lemma \ref{computation-w}, $mw_P=\ord(x_1^\alpha)<m$.
Whence $\ord(\psi_3)=a_3+mw_P=2m$.
But then $\psi_3=\psi_1$, a contradiction.

\subsubsection{}
Thus $i_P=3$, hence $w_P<1$. It gives us that $mw_P=\ord(M)<m$ for some
monomial $M$ with $\wt(M)\equiv 1\bmod m$. Hence
$mw_P=\ord(M)=m-1$ and $w_P=(m-1)/m$.
Futher we may assume that $a_1<a_2$ and $a_1<m$, because $a_3>m$.
Therefore $a_1=3$, $a_2\ge m-3$. For $M$ we have the only possibility
\subsection{}
Thus   $i_P=3$, hence $w_P<1$. It gives us $w_P=(m-1)/m$
and for some monomial $\psi$ we have $\ord(\psi)=m-1$.
Since $a_1=3$, $a_2\ge m-3$ and $a_3\ge m+1$,
$\psi$ can be only
$x_1^k$, where $3k=m-1$.
Whence $m\equiv 1\bmod 3$.
Finally $(X,P)$ is a point of type $cA/m$ because
$m\ge 4$ and $2\not\in\ord(C\3)$.
Thus \ref{t1} and \ref{t2} are proved.

\subsubsection{}
To prove \ref{t3} we remark that $2m+1\in\ord(C\3)$, because
$2m+1\equiv 0\bmod a_1$.
By normalizedness, it gives us $a_3=m+1$ or $a_3=2m+1$.
If $a_3=m+1$, then $2m-3=\frac{m-4}{3}a_1+a_3\in\ord(C\3)$,
whence $a_2=2m-3$ or $m-3$ and we have cases \ref{1)}, \ref{2)}..
So we assume that $a_3=2m+1$. Then  similarly
$3m-3=(m-1)a_1\in\ord(C\3)$, whence $a_2=3m-3$, $2m-3$ or $m-3$.
We get \ref{3)}, \ref{4)} or \ref{5)}.
\par
If in \ref{3)} $m=7$, then  $a_3-7=8\in\ord(C\3)$,
this contradicts normalizedness of the coordinate system
$(x_1,\dots,x_4)$.
In cases \ref{1)}, \ref{3)} $\wt(x_2)\not\equiv\wt(x_3)\bmod m$.
If $(X,P)$ is not a cyclic quotient, then an
equation $\phi=0$ of $X\3$ must contain the term $x_1x_2$
(see \cite{Mori-term}) and therefore another monomial of
order $\ord(x_1x_2)=m$.
Since $a_3>m$, it can be only $x_4$. Therefore $X\3$ is non-singular.
In case \ref{5)} one can obtain representation of $(X,P)$
by using that $C\3$ is a complete intersection in $\CC^4$.
\end{proof}
Now proof of Theorem \ref{pmain} follows from Lemmas
 \ref{uuu} (case \ref{2-cAx/4}),
\ref{cb} (cases \ref{2-IC}, \ref{2-cA/m}) and \ref{cv}
(cases
\ref{3-7}, \ref{3-cyc1}, \ref{3-diff1}, \ref{3-cyc2},
\ref{3-diff2} and \ref{3-diff3}).
Restrictions on the number of Gorenstein singular
points follows from Proposition \ref{grw}.

\end{document}